\documentclass[a4paper]{amsart}
\usepackage{amssymb,amscd,amsthm,verbatim, latexsym}
\usepackage[english]{babel}
\newtheorem{theorem}{Theorem}[section]         
             
\newtheorem{corollary}[theorem]{Corollary}     
\newtheorem{proposition}[theorem]{Proposition} 
         
\newtheorem{remark}[theorem]{Remark}           
\newtheorem{definition}[theorem]{Definition}   
\newtheorem{example}[theorem]{Example}          
\numberwithin{equation}{section}

\newenvironment{namelist}[1]{%
\begin{list}{}
    {
      
      \settowidth{\labelwidth}{#1}
      \setlength{\leftmargin}{1.1\labelwidth}
    }
  }{%
\end{list}}
\newcommand{\nc}{\newcommand}
\nc{\R}{\ensuremath{\mathbb{R}}}
\nc{\K}{\ensuremath{\mathbb{K}}}
\nc{\Z}{\ensuremath{\mathbb{Z}}}
\nc{\Q}{\ensuremath{\mathbb{Q}}}
\nc{\Mo}{\ensuremath{\mathbb{M}}}
\nc{\C}{\ensuremath{\mathfrak {C}}}
\nc{\os}{\ensuremath{\mathfrak {S}}}
\nc{\A}{\ensuremath{\mathcal {A}}}
\nc{\Af}{\ensuremath{\mathfrak {A}}}
\nc{\B}{\ensuremath{\mathfrak {B}}}
\nc{\G}{\ensuremath{\mathfrak {G}}}
\nc{\IND}{\operatorname{IND}}
\nc{\NBC}{\operatorname{NBC}}
\nc{\Ker}{\operatorname{Ker}}
\nc{\OS}{\operatorname{OS}}
\nc{\OT}{\operatorname{OT}}
\nc{\M}{\ensuremath{\mathcal M}}
\nc{\cl}{\ensuremath{c\ell}}
\nc{\sgn}{\operatorname{sgn}}
\nc{\id}{\operatorname{id}}
\begin{document}
\title[Orlik-Solomon type algebras]
{Gr\"obner  and diagonal bases in\\ Orlik-Solomon type algebras}
\thanks{ 
\emph{Secondary}:  14F40.
}
\author{Raul Cordovil and David Forge}
\address{\newline 
Departamento de Matem\'atica,\newline 
Instituto Superior T\' ecnico \newline
 Av.~Rovisco Pais
 - 1049-001 Lisboa  - Portugal}
\email{cordovil@math.ist.utl.pt}
\thanks{The  first  author's research was 
supported in part by FCT (Portugal) through program POCTI and
 the project SAPIENS/36563/00.}
\address{{}\newline
Laboratoire de Recherche en Informatique\newline 
Batiment 490
Universite Paris Sud\newline 
91405 Orsay Cedex - France}
\email{forge@lri.fr}
\begin{abstract} 
The Orlik-Solomon algebra of a matroid $\M$ is the quotient of the 
exterior algebra on the points by the ideal $\Im(\M)$ generated by 
the boundaries of the circuits of the matroid.
There is an isomorphism between the Orlik-Solomon algebra of a complex 
matroid and the cohomology of the complement of a complex arrangement of  
hyperplanes.
In this article a
generalization of the Orlik-Solomon algebras, called $\chi$-algebras, are 
considered. These new algebras  include, apart from the Orlik-Solomon algebras, 
the Orlik-Solomon-Terao algebra of a set of vectors  and the Cordovil 
algebra of an oriented matroid. To encode an important property of the 
``no broken circuit bases"
of the Orlik-Solomon-Terao algebras, Andr\'as Szenes has introduced a 
particular type
of bases, the so called
 ``diagonal bases". This notion
 extends naturally to the $\chi$-algebras.  We give a survey of the 
results  obtained by the  authors concerning the construction 
 of Gr\"obner bases of $\Im_{\chi}(\M)$ and diagonal bases of
Orlik-Solomon type algebras and we present the combinatorial analogue of
an ``iterative residue formula'' introduced by Szenes.
 \end{abstract}
\thanks{Typeset by \AmS -\LaTeX}
\maketitle 
\today
\section{introduction}
Let $\M=\M([n])$ be a matroid on the ground set 
$[n]:=\{1,2,\ldots,n\}$. 
The Orlik-Solomon algebra of a matroid $\M$ is the quotient of the 
exterior algebra on the points by the ideal $\Im(\M)$ generated by 
the boundaries of the circuits of $\M$.
The isomorphism between the Orlik-Solomon algebra of complex matroid and 
the cohomology of the complement of a complex arrangement of  hyperplanes 
was established in 
\cite{OS}. The Orlik-Solomon 
algebras  have been then
intensively studied.
A general reference on hyperplane arrangements and Orlik-Solomon algebras is
\cite{OT}.  Descriptions of developments from the 
early 1980's to the end of 1999, together with the 
contributions of many authors,
can be found in \cite{Falk1,Yuz}.\par
In this article a
generalization of the Orlik-Solomon algebras, called $\chi$-algebra, is 
considered. These new algebras  include, apart from the Orlik-Solomon algebras, 
the Orlik-Solomon-Terao algebra of a set of vectors \cite{OT1} and the 
Cordovil algebra of an oriented matroid \cite{cor2}.
We will survey recent results concerning
this family of Orlik-Solomon type algebras (see \cite{Diagonal,f,FL}). 
In this introduction, we will recall the origin of the Orlik-Solomon
algebra and we will develop the different notions
used in the next sections like matroids and oriented matroids, 
the Orlik-Solomon algebra and its generalizations, its diagonal bases and
the Gr\"obner bases of the defining ideal.\par
Let $V$ be a vector space of dimension $d$  over some field $\mathbb 
{K}$. A (central) arrangement (of hyperplanes) in $V,$
$\mathcal{A}_{\mathbb {K}}=\{H_{1},\ldots,H_{n}\},$ is a
finite listed set of codi\-mension one vector subspaces. 
Given an arrangement $\mathcal{A}_{\mathbb {K}}$\, we always  suppose 
fixed a  family of  linear forms 
$\big\{\theta_{{H}_i}\in V^{*}: H_i\in \mathcal{A}_{\mathbb
{K}},\,\mathrm{Ker}(\theta_{{H}_i})={H_i}
\big\},$ where\, $V^{*}$ denotes the dual space of\,
$V$.  Let 
$L(\A_{\mathbb {K}})$ be  the {\em intersection lattice of
$\A_{\mathbb {K}}$}: i.e., the set of  intersections of 
hyperplanes in $\mathcal{A}_{\mathbb {K}},$
partially ordered by reverse inclusion.
There is a    matroid $\M(\A_{\mathbb {K}})$ on the ground set 
$[n]$  determined by   $\A_{\K}$: a subset $D\subseteq  [n]$\, 
is  a \emph{dependent set} of  $\M(\A_{\mathbb {K}})$ 
iff
there are  scalars $\zeta_{i}\in \mathbb 
{K},\, i\in D$, not all nulls,
such that $\sum_{i\in D}\zeta_{i}\theta_{H_{i}}=0$. A \emph{circuit}  is a 
minimal dependent set with respect to inclusion. 
If $\K$ is an ordered field an additional 
structure is obtained: to every circuit $C,$ 
$\sum_{i\in C}\zeta_{i}\theta_{H_{i}}=0,$ we associate a
partition (determined up to a factor $\pm 1)$  
$C^+:=\{i\in C: \zeta_{i}>0\}, C^-:=\{i\in C:
\zeta_{i}<0\}$. With this new structure $\M(\A_{\mathbb {K}})$ is said a 
\emph{(realizable) oriented matroid} and denoted by
$\boldsymbol{\M}(\A_{\mathbb {K}})$.  Oriented
matroids on a ground set $[n]$, denoted $\boldsymbol{\mathcal{M}}([n]),$ 
are a very
natural mathematical concept and can be seen as the theory of
generalized hyperplane arrangements, see
\cite{MO}.\par
Set $\mathfrak{M}(\mathcal{A}_{\mathbb{K}}):=
V\setminus\bigcup_{H\in
\mathcal{A}_{\mathbb{K}}} H$. 
The manifold \,$\mathfrak{M}(\mathcal{A}_{\mathbb{C}})$\, plays an 
important
role in the 
Aomoto-Gelfand theory
of multidimensional hypergeometric functions (see \cite{OT2} for 
a recent introduction
from the point of  view of arrangement theory).  Let ${K}$ be a 
commutative ring. In 
\cite{OS,OS1,OT}  the determination of the cohomology $K$-algebra
$H^{*}\big(\mathfrak{M}(\mathcal{A}_{\mathbb{C}}); K\big)$\, from the 
matroid
$\mathcal{M}(\mathcal{A}_{\mathbb {C}})$ is accomplished by first 
defining the Orlik-Solomon $K$-algebra 
$\mbox{OS}(\mathcal{A}_{\mathbb{C}})$
in terms of generators and 
relators which depends only on the matroid $\mathcal{M}
(\mathcal{A}_{\mathbb{C}})$, and then by
showing that this
algebra is isomorphic to
$H^{*}\big(\mathfrak{M}(\mathcal{A}_{\mathbb{C}}); K\big)$.
Aomoto suggested the study of the (graded) $\mathbb{K}$-vector space 
$\mbox{AO}
({\mathcal{A}}_{\mathbb{K}}),$ 
generated by the basis 
$\{Q({\mathcal{B}}_{\mbox{I}})^{-1}\},$ where $I$  is an independent
set of  
$\mathcal{M}
({\mathcal{A}}_{\mathbb{K}}),$   ${\mathcal{B}}_{\mbox {I}}:= \{H_i\in
\mathcal{A}_{\mathbb{K}}:i\in I\},$\,   and 
\,$Q ({\mathcal{B}}_{\mbox{I}})=\prod_{i\in I}\theta_{H_i}$ 
denotes the corresponding defining polynomial. To 
answer a conjecture of Aomoto,
Orlik and Terao have  introduced in \cite{OT1} a commutative 
$\mathbb{K}$-algebra, $\OT({\mathcal{A}}_{\mathbb{K}}),$ called the Orlik-Solomon-Terao
algebra. The algebra $\OT({\mathcal{A}}_{\mathbb{K}})$ is isomorphic to 
$\mbox{AO}(\mathcal{A}_{\mathbb {K}})$
 as a graded  $\mathbb {K}$-vector space in terms of the equations 
$\{\theta_{H}: H\in \mathcal{A}_{\mathbb{K}}\}$.
A  ``combinatorial analogue'' of the   
algebra of Orlik-Solomon-Terao  was introduced in \cite{cor2}:
to every oriented matroid $\boldsymbol{\mathcal{M}}$ was associated  
a commutative
$\mathbb{Z}$-algebra, denoted by
$\mathbb{A}(\boldsymbol{\mathcal{M}})$ and called the Cordovil algebra. 
The $\chi$-algebras generalizes
the three just mentioned algebras: Orlik-Solomon, Orlik-Solomon-Terao 
and the Cordovil algebras, see \cite{FL} or 
Example~\ref{chi} below.\par
In section two we will give the definition of a $\chi$-algebra 
and recall the principal examples. In general a $\chi$-algebra, denoted 
$\mathbb{A}_\chi({\mathcal{M}}),$ is defined as the quotient of
some kind of a finite $\K$-algebra $\Af$ by an ideal $\Im_{\chi}(\M)$ of
$\Af$   whose
generators  are defined from the circuits of $\M$ and are depending of the map $\chi$, see
Definition~\ref{FL1}. In particular  the first important result is that like for the original
Orlik-Solomon algebra we get $\boldsymbol{nbc}$-bases of the $\chi$-algebra (as a module) from the 
``no broken circuit''
sets of the matroid and corresponding basis for the ideal 
$\Im_{\chi}(\M)$.\par
In section three, we construct the reduced Gr\"obner basis of the ideal 
$\Im_{\chi}(\M)$ for any term order $\prec$ on the set of the monomials
$\mathbb{T}(\Af)$ of the algebra $\Af$. This result gives as a corollary a universal  Gr\"obner
basis (a Gr\"obner basis who works for every term order)  which is shown
to be minimal. Finally we remark  that the $\boldsymbol{nbc}$-bases are in some sense
the bases corresponding to the  Gr\"obner bases for the different term
orders.\par
In section four, following Szenes \cite{Sz}, we define a particular type
of basis of $\mathbb{A}_\chi,$ the so called ``diagonal basis", see
Definition
\ref{diagonal}. The $\boldsymbol{nbc}$-bases are an important examples of diagonal bases. We construct
the dual bases  of these bases, see Theorem~\ref{dbases}.
 Our definitions make also use of an ``iterative
residue formula" based on the matroidal operation of contraction, see 
Equation
(\ref{resid}).  This formula can be seen as the combinatorial 
analogue of an ``iterative residue formula'' introduced by Szenes, 
\cite{Sz}.
As applications we deduce nice formulas to express
a pure element in a diagonal basis.
We prove also that the $\chi$-algebras verify a splitting short  exact
sequence, see Theorem~\ref{thm: There is}.
This theorem generalizes for the $\chi$-algebras previous similar theorems 
of
 \cite{cor2, OT}. \par
We use \cite{W1, W2} as a general reference in matroid theory. We refer to 
\cite{MO} and \cite{OT} for good sources of the theory
of oriented matroids and arrangements of hyperplanes,  respectively.
\section{$\chi$- algebras}
Let $\mbox{IND}_\ell(\mathcal{M})\subseteq\binom{[n]}{\ell}$ 
[resp. $\mathrm{DEP}_\ell(\mathcal{M})\subseteq\binom{[n]}{\ell}$]
be the 
family of 
independent [resp. dependents] sets  of cardinality $\ell$ of the matroid $\mathcal{M}$
and    set\,
$$\mbox{IND}(\mathcal{M}) :=\bigcup_{\ell \in 
\mathbb{N}}\mathrm{IND}_\ell(\mathcal{M}),$$
$$\mbox{DEP}(\mathcal{M}) :=\bigcup_{\ell \in 
\mathbb{N}}\mathrm{DEP}_\ell(\mathcal{M}).$$
 We denote by $\mathfrak{C}=\mathfrak{C}(\mathcal{M})$ the set of 
circuits  of $\mathcal{M}$. For shortening of the notation
the singleton set  $\{x\}$ is denoted  by  $x$.
When the smallest element $\alpha$ of a circuit $C,$ $|C|>1,$ is deleted, 
the
remaining  set,  $C\setminus \alpha,$\, is said to be a
\emph{broken circuit}. 
(Note that our definition is slightly different from the 
standard one. In the standard definition
$C\setminus \alpha$ can be empty.) 
A \emph{no broken circuit} set of a matroid $\mathcal{M}$ 
is an independent subset of $[n]$ which does not contain any broken 
circuit. 
Let $\mathrm{NBC}_{\ell} (\mathcal{M})\subseteq\binom{[n]}
{\ell}$ be the set of the no broken
circuit sets of cardinal $\ell$ of $\mathcal{M}$ and set   
$$\mbox{NBC}(\mathcal{M})
:=\bigcup_{\ell \in \mathbb{N}}\mbox{NBC}_\ell(\mathcal{M}).$$  Let
 $L(\mathcal{M})$ be the lattice of  flats of
$\mathcal{M}$.  \big(We remark that the 
lattice map\,  $\phi: L(\mathcal{A}_{\mathbb {K}})\to L(\mathcal{M}
(\mathcal{A}_{\mathbb{K}}))$, determined by  the one-to-one correspondence 
$\phi':H_{i}\longleftrightarrow \{i\}$, $i=1,\dotsc, n,$
is a lattice isomorphism.\big{)}  For an   independent set $I\in \mathrm{IND}(\mathcal{M})$,  
let  $\cl(I)$ be the closure of $I$ in $\mathcal{M}$.\par
For every permutation $\sigma \in$
${\os}_m,$  let $X^\sigma$ be the ordered set $$X^\sigma:=i_{\sigma(1)}\prec\cdots \prec
i_{\sigma(m)} 
=(i_{\sigma(1)},i_{\sigma(2)},\,\ldots,\, i_{\sigma(m)}).$$  
 When necessary we also see the set $X=\{i_1,\dotsc, i_m\},$ as the ordered set
$$X^{\mathrm{id}}=(i_1,\dotsc, i_m).$$ Set $X^\sigma\setminus
x:=(i_{\sigma(1)},\dotsc,\widehat{x},\dotsc,i_{\sigma(m)}
)$. If\, $Y^\beta=(j_{\beta(1)},\dotsc, j_{\beta(m')})$  and 
$X\cap Y=\emptyset,$ let $X^\sigma\circ Y^\beta$ be the concatenation
of   $X^\sigma$ and $Y^\beta$, i.e., the ordered set
$$
X^\sigma\circ Y^\beta:=
(i_{\sigma(1)},\dotsc,
i_{\sigma(m)},j_{\beta(1)},\dotsc, j_{\beta(m')}).$$ 
\begin{definition}\label{algebra}
{\em
Let $\chi$ be a mapping $\chi: 2^{[n]}\to \K$.  
Let us also define $\chi$ for 
ordered sets by $\chi (X^\sigma )= \sgn(\sigma) \chi (X)$, 
where  $\sgn(\sigma)$ denotes 
the sign of the permutation $\sigma$.
Fix a set $E=\{e_1,\ldots, e_n\}$.
Let  $\Af:=\K\oplus\Af_1\oplus\cdots\oplus \Af_n$ be the graded 
algebra over the field\, $\K$ generated by the
elements $1, e_1,\dotsc,e_n$  and satisfying the relations:
\begin{namelist}{xxxx}
\item[~$\circ$]$ 1e_i=e_i1=e_i$, for all\, $e_i\in E$,
\item[~$\circ$] $e_i^2=0$, for all\, $e_i\in E$,
\item[~$\circ$] $e_j\cdot e_i=\beta_{i,j}e_i\cdot e_j$\, with\,
$\beta_{i,j}\in 
\K^*$\,\, for all\,
$i<j$.
\end{namelist}
By definition the $\chi$-{\it boundary}
of an element $e_X\in \Af$, $X\not =\emptyset$,  is   given by the formula
$$\partial e_X:=\sum_{p=1}^{p=m}(-1)^p
\chi(X\setminus i_p)e_{X\setminus i_p}.$$
We set $\partial e_i=1$,
for all
$e_i \in E$.
We extend $\partial$ to the $\K$-algebra $\Af$ by linearity. 
}
\end{definition}
Let $X=(i_{1}, i_{2},\dotsc, i_{m})$. In the 
sequel we will denote by
$e_{X}$ the   (pure) element of the $\K$-algebra \, $\Af$,
$$e_{X}:=e_{i_1}\cdot e_{i_2}\cdot \cdots \cdot e_{i_m}.$$
By
convention we
 set  $e_{\emptyset}:=1$. 
Both the  exterior $\K$-algebra,\, $\bigwedge E$,\, (take
$\beta_{i,j}=-1$) and the  polynomial algebra $\K[e_{1},\dotsc, e_{n}]/\langle e_{i}^{2}\rangle$ with squares zero 
(take  $\beta_{i,j}=1$) considered in \cite{cor2, OT1}, are such $\K$-algebras $\Af$
and will be the only ones to be used in the examples.
It is clear that for any $x\not \in X,$ $$\pm \partial e_{X\cup x}=
(-1)^{m+1}\chi (X)e_{X}+\sum_{p=1}^{p=m}(-1)^p
\chi(X\setminus i_{p}\circ x)e_{X\setminus i_{p}\cup x}.$$
From the equality $\chi (X^\sigma )= \sgn(\sigma) \chi (X)$,
it is easy to see that for $\sigma \in \os_{|X|}$ we have $$\partial e_X=
\sgn(\sigma)\sum_{p=1}^{p=m}(-1)^p
\chi(X^\sigma \setminus i_{\sigma (p)})e_{X\setminus i_{\sigma (p)}}.$$
Given an independent set $I,$ an element $a\in \cl(I)\setminus I$
is said 
\emph{active} in $I$ if $a$ is the minimal element of the 
unique circuit contained in
$I\cup a$. We say that a subset $U\subseteq [n]$ is a 
\emph{unidependent} set of  $\M([n])$ if it
contains a unique circuit, denoted $C(U)$. 
 Note that
$U$ is unidependent iff
$\mathrm{rk}(U)=|U|-1$. We say that a unidependent set 
$U$ is an \emph{inactive
unidependent} if $\mathrm{min} (C(U))$ is the the minimal active element 
of 
$U \setminus\mathrm{min} (C(U))$. We will denote by 
$\mathrm{UNI}_\ell(\mathcal{M})$ 
for the sets of inactive unidependent sets of size $\ell$ and  set
$$\mbox{UNI}(\mathcal{M})
:=\bigcup_{\ell \in \mathbb{N}}\mbox{UNI}_\ell(\mathcal{M}).$$
Let us
remark that   $U$ is a unidependent set of  $\M$ iff for some (or every) $x\in
U,$
$\mathrm{rk}(x)\not =0,$
$U\setminus x$ is a unidependent set of $\M/x$. 
\begin{definition}[\cite{FL}]\label{FL1}
{\em Let $\chi$ be a mapping $\chi: 2^{[n]}\to \K$.
Let $\Im_{\chi}(\M([n]))$ be the (right) ideal of\, $\Af$
generated by the
$\chi$-boundaries $\{\partial e_C: C\in \C(\M), |C|>1\}$ and the set
of the loops of $\M$, 
$\{e_i: \{i\}\in
\C(\M)\}$.  We say that
$\mathbb{A}_\chi(\M):=\Af/{\Im}_{\chi}(\M)$
is a $\chi$-{\it algebra} if $\chi$ satisfies the
following two properties:
\begin{namelist}{xxxxxxxx}
\item[~$(\ref{FL1}.1)$]
$\chi(I)\not=0$ if and only if $I$ is independent.
\item[~$(\ref{FL1}.2)$]
For any two unidependents $U$  and $U'$ of $\M$ with $U'\subseteq U$ 
there is a  scalar $\varepsilon_{_{U,U'}}\in \K^*,$  such that
$\partial e_U=\varepsilon_{_{U,U'}}(\partial e_{U'})e_{U\setminus U'}$.
\end{namelist}
}
\end{definition}
Note that
$$
\{e_C:\, C\in  \C(\M)\} \subseteq \Im_{\chi}(\M([n])).
$$
For every $X\subseteq [n],$ we denote by $[X]_{\mathbb{A}}$ or
shortly by   $e_X$ when no confusion will result, the residue 
class in $\mathbb{A}_\chi({\mathcal{M}})$ 
determined by the  element $e_X$.  
Since ${\Im}_{\chi}(\M)$ is a homogeneous  ideal,   
$\mathbb{A}_\chi(\M)$ inherits a grading from $\Af$. More precisely
we have
$\mathbb{A}_\chi({\mathcal{M}})=
\K\oplus\mathbb{A
}_1\oplus
\cdots 
\oplus 
\mathbb{A}_r,$ 
where\, $\mathbb{A}_\ell=\Af_\ell/\Af_{\ell}\cap{\Im}_{\chi}(\M)$ denotes 
the 
subspace  of $\mathbb{A}_\chi({\mathcal{M}})$ generated by the elements 
$\big\{[I]_{\mathbb{A}}: I\in 
\mbox{IND}_\ell({\mathcal{M}})\big\}$.
Set 
\begin{equation*}
\boldsymbol{nbc}_{\,\ell}:=\big\{[I]_{\mathbb{A}}:  
I\in \mbox{NBC}_{\,\ell}({\mathcal{M}})\big\} \quad\mbox{and}\quad
\boldsymbol{nbc}:=\bigcup_{\ell=0}\boldsymbol{nbc}_{\,\ell},
\end{equation*}
\begin{equation*}
\boldsymbol{dep}_{\ell}:=\big\{[D]_{\mathbb{A}}:  
D\in \mbox{DEP}_{\ell}({\mathcal{M}})\big\} \quad\mbox{and}\quad
\boldsymbol{dep}:=\bigcup_{\ell=0}\boldsymbol{dep}_{\,\ell},
\end{equation*}
\begin{equation*}
\boldsymbol{uni}_{\ell}:=\big\{[U]_{\mathbb{A}}:  
U\in \mbox{UNI}_{\ell}({\mathcal{M}})\big\}\, \quad\mbox{and}\quad
\boldsymbol{uni}:=\bigcup_{\ell=0}\boldsymbol{uni}_{\,\ell}.
\end{equation*}
\begin{remark}\label{rem1}
{\em
From $(\ref{FL1}.1)$ and $(\ref{FL1}.2)$ we conclude that
${\Im}_{\chi}(\M)$ has the basis
$\boldsymbol{dep}\cup\partial \boldsymbol{uni}$
and that 
$\boldsymbol{nbc}:=\big\{[I]_{\mathbb{A}}:   I\in
\mbox{NBC}({\mathcal{M}})\big\}$ is a basis of 
 the vector space $\mathbb{A}=\mathbb{A}_\chi(\M)$. 
We also have that  $\boldsymbol{nbc}_\ell$ is
a   basis of the vector space $\mathbb{A}_\ell$. This 
fundamental property was first discovered for
the Orlik-Solomon algebras \cite{OT}, and
then also for the other  classical $\chi$-algebras, see \cite{cor2,
OT1}  and the following example for more details. Note also that this 
implies
that
$[X]_{\mathbb{A}}\not=0$ iff $X$  is an independent set of
${\mathcal{M}}$.
}
\end{remark}
\begin{example}\label{chi}
{\em
Recall the three usual $\chi$-algebras $\mathbb{A}_\chi(\M)$.
\begin{namelist}{xxxx}
\item[~~$\circ$] Let $\Af=\bigwedge E$ be the  exterior $\K$-algebra 
(taking $\beta_{i,j}=-1$).  Setting
$\chi(I^\sigma)=\mbox{sgn}(\sigma)$   for every independent set 
$I$ of a matroid $\M$ and
every permutation
$\sigma \in \os_{|I|},$ we obtain the Orlik-Solomon algebra,
$\OS(\M)$.
\item[~~$\circ$]
Let $\A_\K=\{H_i: H_i=\Ker(\theta_i), i=1,2,\dotsc, n\}$ be 
an hyperplane arrangement
and
$\M(\A_\K)$  its associated  matroid.
For every flat $F:=\{f_1,\dotsc, f_k\}\subseteq [n]$ of $\M(\A_\K)$ 
we choose a bases $B_F$
of the vector subspace of $(\K^d)^*$ generated by 
$\{\theta_{f_1},\dotsc, \theta_{f_k}\}$. By taking  
$\Af=\K[e_{1},\dotsc, e_{n}]/\langle e_{i}^{2}\rangle$ the  polynomial
algebra with squares null  (taking $\beta_{i,j}=1$) and taking for any
$\{i_1,\dotsc,i_\ell\}=I\in
\mathrm{IND}_\ell,$\,\,
$\chi(I)=\det(\theta_{i_1},\dotsc,\theta_{i_\ell}),$ where the 
vectors are expressed in
the
  basis $B_{\cl(I)},$ we obtain the Orlik-Solomon-Terao algebra $\OT(\A_\K),$
defined in
\cite{OT1}.
\item[~~$\circ$]
Let $\boldsymbol{\M}([n])$ be an oriented matroid. For every 
flat $F$ of $\boldsymbol{\M}([n]),$ we choose
(determined up to a factor $\pm 1$) a bases signature in the 
restriction of $\boldsymbol{\M}([n])$ to $F$.
We define a \emph{signature of the independents of an oriented 
matroid $\boldsymbol{\M}([n])$} as a mapping,
${\sgn}:\mbox{IND}(\boldsymbol{\M})\to \{\pm 1\},$ where 
$\sgn(I)$ is equal to the basis signature of $I$ in
the restriction of $\boldsymbol{\M}([n])$ to $\cl(I)$.
By taking 
 $\Af=\Q[e_{1},\dotsc, e_{n}]/\langle e_{i}^{2}\rangle$ the  polynomial algebra over the 
rational field $\Q$ with squares zero 
 (take  $\beta_{i,j}=1$) and taking $\chi(I)=\sgn(I)$ (resp.  
$\chi(X)=0$)
 for every independent (resp. dependent) set of the matroid, 
we obtain the algebra
$\mathbb{A}(\boldsymbol{\M})\oplus_\Z \Q,$ where
$\mathbb{A}(\boldsymbol{\M})$ denotes the Cordovil $\Z$-algebra defined in
\cite{cor2}.
\end{namelist}
}
\end{example}
\section{Gr\"obner bases of $\chi$-ideals}
For general details on Gr\"obner bases of an ideal, see \cite{Adams,Becker}. We begin by
adapting some  definitions to our context. 
Consider the $\K$-algebra $\Af$ introduced in Definition~\ref{algebra}.
Note that there are monomials $e_Y, e_Z\in
\Af$, such that
$e_Y\cdot e_Z=0$. In the standard case where
 $\Af$ is replaced by the polynomial ring $\K[e_{1},\dotsc, e_{n}]$, 
this is not
possible. So the the following definitions are slightly different
from the standard corresponding ones given in \cite{Adams,Becker}. Let 
$\M=\M([n])$ be a matroid,
${\Im}_{\chi}(\M)$ and 
${\mathbb A}_{\chi}(\M)$
the $\chi$-ideal and $\chi$-algebra as defined in the previous section. 
We will denote for shortness 
${\mathbb A}(\M)$ for 
${\mathbb A}_{\chi}(\M)$.
\par
\begin{definition}
{\em
Let $\mathbb{T}=\mathbb{T}(\Af)$ be the \emph{set of the 
monomials of the
$\K$-algebra $\Af$}, i.e.,
$\mathbb{T}(\Af):=\{e_X:\, X=(e_{i_1},\dotsc,e_{i_m})\}$.
A total ordering $\prec$\, on the monomials $\mathbb{T}$   
  is said a
\emph{term order} on $\mathbb{T}$ if
$e_\emptyset=1$ is the minimal element and $\prec$  is compatible with the multiplication in $\Af$, i.e.,   $$\forall 
e_X,e_Y,e_Z\in \mathbb{T},\hspace{2mm} (e_X\prec e_Y)\&
(e_X\cdot e_Z\not =0)\& (e_Y\cdot e_Z\not =0)\, \Longrightarrow\, 
e_{X}\cdot e_{Z}\prec e_{Y}\cdot e_{Z}.$$
Given a term order\, $\prec$\, on $\mathbb{T}$ and a non-null polynomial $f\in 
\Af$, 
we may write
$$f=a_1e_{X_1}+a_2e_{X_2}+\dotsm+a_me_{X_m},$$
where $a_i\in \K^{*}$ and $e_{X_m}\prec\dotsm\prec e_{X_1}$. 
We say that the
$a_ie_{X_i}$ [resp. $e_{X_i}$] are the \emph{terms} [resp.
\emph{monomials}] of $f$. We say that
$\mbox{lp}_\prec(f):=e_{X_1}$ [resp.
$\mbox{lt}_\prec(f):=a_1e_{X_1}$] is the \emph{leading monomial } 
[resp. \emph{leading term}] of $f$ (with respect to $\prec)$.
We also define $\mbox{lp}_\prec(0)=\mbox{lt}_\prec(0)=0$. 
Note that in general we have $\mbox{
lp}_\prec(hg)\not=\mbox{lp}_\prec(h)\mbox{lp}_\prec(g),$ 
contrarily to the cases considered in \cite{Adams,Becker}.}
\end{definition}
\begin{example}
{\em 
A permutation $\pi\in \Sigma_n$ 
defines a linear reordering of the set $[n]$:
$\pi^{-1}(1)<_\pi\pi^{-1}(2)<_\pi\dotsm <_\pi\pi^{-1}(n)$. Consider the
ordering of the set $E$
$$e_{\pi^{-1}(1)}\prec_\pi\, e_{\pi^{-1}(2)}\prec_\pi\dotsm \prec_\pi
e_{\pi^{-1}(n)}.$$  The corresponding
degree lexicographic ordering on the monomials
$\mathbb{T}$, also denoted $\prec_\pi$, is a
term order on $\mathbb{T}$.
 }
\end{example}
For a subset $S,\, S\subseteq \Af$ and a term order\, $\prec$\, 
on $\mathbb{T}(\Af),$ we
define the 
\emph{leading term
ideal} of $S$, denoted $\mbox{Lt}_\prec(S)$, as the 
ideal generated by the leading
monomials of the polynomial in
$S$, i.e.,
$$\mbox{Lt}_\prec(S):=\langle \mbox{lp}_\prec(f):\, f\in S\rangle.$$
 In the remaining of this section we suppose that $\M([n])$ is a 
loop free matroid.
\begin{definition}\label{ideal} 
{\em Let $\M$ be a matroid. Let $\prec$ be a 
term order on $\mathbb{T}(\Af)$. Consider the ideal
$\Im_{\chi}(\M)$  of $\Af$ 
A  family  $\mathcal{G}$ of non-null polynomials of the ideal 
$\Im_{\chi}(\M)$ is called
a 
\emph{Gr\"obner basis} of the ideal $\Im_{\chi}(\M)$ with respect to
$\prec$\, iff    $$\mbox{
Lt}_\prec(\mathcal{G})=\mbox{Lt}_\prec(\Im_{\chi}(\M)).$$  The Gr\"obner
basis $\mathcal{G}$ is called  \emph{reduced}   
 if, for every element\, $g\in \mathcal{G}$ we have 
$\mbox{lt}_\prec(g)=\mbox{lp}_\prec(g),$ and for every two distinct
elements $g,g'\in \mathcal{G},$  no term of
$g'$ is divisible by $\mbox{lp}_\prec(g)$. The Gr\"obner basis
$\mathcal{G}$ is called a
\emph{universal Gr\"obner basis} if it is  a Gr\"obner basis  with respect to all term
orders on $\mathbb{T}(\Af)$  simultaneously. If $\mathcal{U}$ is a universal Gr\"obner basis, minimal for
inclusion with this property, we say that  $\mathcal{U}$ is a \emph{minimal universal 
Gr\"obner basis}.}
\end{definition}
From Definition~\ref{ideal} we conclude:
\begin{proposition} \label{David} Let\, $\mathcal{G}_{\prec}$ be a  
Gr\"obner basis   of\, the  ideal
$\Im_{\chi}(\M)$ with respect to the term order $\prec$ on\, $\mathbb{T}(\Af)$.
Then 
$$\mathcal{B}_{\mathcal{G}_{\prec}}:=\big\{e_X+\Im_{\chi}(\M):\, 
e_X\not \in 
\mbox{\em Lt}_\prec(\mathcal{G})\big\}$$ 
is a  basis of the module
 $\mathbb{A}_\chi (\M)$. \qed
\end{proposition}
We say that the well determined basis\, $\mathcal{B}_{\mathcal{G}_{\prec}}$ is the
{\em canonical basis  of the $\chi$-algebra
$\mathbb{A}_\chi (\M)$ for
the Gr\"obner basis $\mathcal{G}$ of the  ideal $\Im_{\chi}(\M)$,
with respect to the term order\, $\prec$\, on\, $\mathbb{T}(\Af)$.}\par
Consider the partition $\mathbb{T}(\Af)=\mathbb{T}_{i}(\Af)\biguplus \mathbb{T}_{d}(\Af)$
where: 
$$\mathbb{T}_{i}(\Af):=
\big\{e_I:\, I\in \mbox{IND}(\M)\big\}\,\, \,\,\mbox{and}\,\,\,\,
\mathbb{T}_{d}(\Af):=\big\{e_D:\, D\,\, \in \mbox{DEP}(\M) \big\}.$$
Let $\K[\mathbb{T}_{i}]$  and  $\K[\mathbb{T}_{d}]$  be the
$\K$-vector subspaces of $\Af$ generated  by the basis 
 $\mathbb{T}_{i}$ and $\mathbb{T}_{d},$ respectively. So \, $\Af=
\K[\mathbb{T}_{i}]\,\oplus\, \K[\mathbb{T}_{d}]$. We also see  the set
$\K[\mathbb{T}_{d}]\subseteq
\Im_{\chi}(\M)$
as the  ideal of
$\Af$ generated by the set of monomials
$\{e_C:\, C\in
\mathfrak{C}(\M)\}$.
 Let $\boldsymbol{p_i}: \Af\to \K[\mathbb{T}_i]$ be the first  projection. 
We define  the term orders on the set of monomials
$\mathbb{T}_{i}$   in a similar way to the  corresponding
definition on  
$\mathbb{T}$. It is clear that the restriction of every term order
 of
$\mathbb{T}$ to the subset $\mathbb{T}_{i}$ is also a term order on
$\mathbb{T}_{i}$.  We can also add to 
$\K[\mathbb{T}_{i}]$  a structure of
$\K$-algebra with the product
$\star : \K[\mathbb{T}_{i}]\times \K[\mathbb{T}_{i}]\to 
\K[\mathbb{T}_{i}],$ 
determined by the equalities
$$e_X\star e_{X^{\prime}}= \boldsymbol{p_i}(e_X e_{X^{\prime}})\quad\mbox{for all}\quad
X, {X^{\prime}}\in
\Im_{\chi}(\M).$$ 
Note that if  $e_X\star e_{X^{\prime}}\not =0$, then $e_X\star
e_{X^{\prime}}=e_X e_{X^{\prime}}$.  We remember that  $e_Xe_{X^{\prime}}\not=0$ iff
$X\cap X^{\prime}=\emptyset$  and $X\cup X^{\prime}\in \mbox{IND}(\M)$. So
$\Im_{\chi_{i}}(\M):= \boldsymbol{p_i}\big(\Im_{\chi}(\M)\big)$ is an ideal of 
$\K[\mathbb{T}_{i}]$.
\begin{proposition}\label{David4}
 Let $\prec$ be a term order on\,\, $\mathbb{T}(\Af)$. Then  the 
leading term
ideals of\, $\Af$, $\mathrm{Lt}_\prec(\boldsymbol{p_i}(\Im_{{\chi}}(\M)))$ and\,
$\mathrm{Lt}_\prec(\Im_{{\chi}}(\M))$ are equal.
In particular a
 Gr\"obner basis of the ideal \, $\Im_{\chi_{i}}(\M)$ of\, $\K[\mathbb{T}_{i}]$ with
respect to term order $\prec$ on $\mathbb{T}_{i}$ is  also a 
Gr\"obner basis of the ideal \, $\Im_{\chi}(\M)$  of\, $\Af$ with respect to the term
order 
$\prec$ on $\mathbb{T}$.
\end{proposition}
\begin{proof}
Note first that if we see $\Im_{\chi}(\M)$ as   a
$\K$-vector space it is clear that 
$\Im_{\chi}(\M)=\Im_{\chi_{i}}(\M)\oplus
\K[\mathbb{T}_{d}]$.
Pick a
non-null polynomial $f\in \Im_{\chi}(\M)$ and let $e_{X_1}:=\mbox{lp}_\prec(f)$. 
 So $e_{X_1}\in \Im_{i}(\M)$ if $X_1\in
\mbox{IND}(\M),$ or\, $e_{X_1}\in \K[\mathbb{T}_{d}]\setminus 0$ if\,
$X_1$ is   a dependent set of $\M$.  If $X_1\in
\mbox{IND}(\M)$  then $e_{X_1}\in \mathrm{Lt}_\prec(\Im_{\chi}(\M))$.
Suppose now that $X_1$ is a dependent set
of $\M$. Then there is a circuit
 $C\subseteq X_1$.  From Definition~\ref{FL1} we know that
$\partial e_C\in \Im_{\chi}(\M)$. It is clear that 
$e_C\in \mathrm{Lt}_\prec(\boldsymbol{p_i}(\Im_{{\chi}}(\M)))$ 
 and  so we have also
$e_{X_1}\in
\mathrm{Lt}_\prec(\boldsymbol{p_i}(\Im_{{\chi}}(\M)))$.
\end{proof}
\begin{remark} \label{David2}{\em  
It is well known that a term order $\prec$ of $\mathbb{T}(\Af)$ 
determines also
a unique reduced Gr\"obner basis of $\Im_{\chi}(\M)$  denoted 
$(\mathcal{G}_r)_{\prec}$. From the definitions we can deduce
also that, for every pair of term orders $\prec$ and $\prec^{\prime}$ on 
$\mathbb{T}(\Af)$,
 $$\mathcal{B}_{\mathcal{G}_{\prec}}= 
\mathcal{B}_{\mathcal{G}_{\prec^{\prime}}}\Leftrightarrow(\mathcal{G}_r)_\prec
=(\mathcal{G}_r)_{\prec^{\prime}}
\Leftrightarrow\mbox{Lt}_\prec\big(\Im_{\chi}(\M)\big)=
\mbox{Lt}_{\prec'}\big(\Im_{\chi}(\M)\big).$$
} 
\end{remark}
\begin{definition}
{\em
For a term order $\prec$ on $\mathbb{T}(\Af)$  
we say that  $\pi_\prec\in \os_n,$ is the 
\emph{ 
 permutation compatible} with  $\prec$ if, for every pair $i,j\in [n],$ 
we have $$e_i\prec e_j\quad \mbox{iff} \quad i<_{\pi_\prec} j \,\,
\big(\Leftrightarrow\,  {\pi_\prec}^{-1}(i)<
{\pi_\prec}^{-1}(j)\,\big).$$  
}
\end{definition}
Let $\mathfrak{C}_{\pi_\prec}$ be the
subset of circuits of $\M$ such that:
\begin{namelist}{xxx}
\item[~$\circ$]  $C\in \mathfrak{C}_{\pi_\prec}$\, iff\,
$\mathrm{inf}_{<_{\pi_\prec}}(C)=
\alpha_\pi(C)\, 
\big(\!=\mathrm{inf}_{<_{\pi_\prec}}
(\mathrm{cl}(C)\setminus C)
\big)$ 
and
$C\setminus\alpha_\pi(C)$ is inclusion minimal   with this property.
\end{namelist} 
In the following we replace  $``\pi_\prec"$ by $``\pi"$ if no mistake can  
results.
We recall that  given a unidependent set $U$ of the matroid $\M([n])$, 
$C(U)$ denotes the unique circuit 
of $\M$ contained in $U$.
\begin{theorem}\label{red}
Let $\prec$ be a term order on\, $\mathbb{T}(\Af)$ compatible with the 
permutation 
$\pi\in \os_n$. Then the family\, 
$\mathcal{G}_{r}:=\big\{\partial e_{C(U)}:\, U\in 
\mathfrak{C}_{\pi_\prec}(\M)\big\}$
form a reduced Gr\"obner
basis  of\, $\Im_{\chi_{i}}(\M)$ with respect to the term order $\prec$.
\end{theorem}
\begin{proof}
From Proposition~\ref{David4} it is enough to prove that 
$(\mathcal{G}_{r})_\prec$ is a reduced Gr\"obner
of $\Im_{\chi_i}(\M)$.
Let $f$ be any element of $\Im_{\chi_{i}}(\M)$,
we have from  Theorem \ref{rem1} that $$f=
\sum _{U\in \mathfrak{U}_{\pi}}
\xi_U \partial e_U, \,\, \xi_U\in \K^\star.$$ Let now remark that
$\mbox{lp}_\prec\big(\partial e_U \big)=e_{U\setminus \alpha_\pi (U)}$
and that these terms are all different.
We have then clearly that $$\mbox{lp}_\prec(f)=\mathrm{sup}_\prec
\big\{\mbox{lp}_\prec(\partial e_{U}):\,\, U\in \mathfrak{U}_{\pi}\big\}.$$
Given an arbitrary $\,U^{\prime}\in 
 \mathfrak{U}_\pi(\M)$ it is clear that $\alpha_{\pi} (C(U^{\prime}))=\alpha_{\pi}
(U^{\prime})$. So,
$$C(U^{\prime})\setminus\alpha_\pi(C(U^{\prime}))\subseteq U^{\prime}\setminus \alpha_{\pi}(U^{\prime}).$$ 
Let $C^{\prime}$ be a circuit of 
$\mathfrak{C}_{\pi}$ such that $C'\setminus \alpha_\pi(C^{\prime})\subseteq 
C(U)\setminus \alpha_\pi(C(U))$.
So we have 
that $\mbox{lp}_\prec(\partial e_{C^{\prime}})$ divides 
$\mbox{lp}_\prec(\partial e_U),$ and $(\mathcal{G}_{r})_\prec$ is a
Gr\"obner basis.\par
Suppose for a contradiction that $(\mathcal{G}_{r})_\prec$ is not a  
reduced Gr\"obner basis: i.e., there exists two
circuits $C$ and $C^{\prime}$ in  
$\mathfrak{C}_{\pi}$ and an element
$c\in C$ such that 
$e_{C^{\prime}\setminus\alpha_\pi(C^{\prime})}$ divides $e_{C\setminus c}$ 
$\big(\Leftrightarrow C'\setminus 
\alpha_\pi(C')\subseteq C\setminus c\big)$.   First we can say 
that $c\not= \alpha_\pi(C)$ because the sets $C^{\prime}\setminus \alpha_\pi(C')$ 
and
$C\setminus \alpha_\pi(C)$
are incomparable. This in particular implies that 
$\alpha_\pi(C)\in C'\setminus \alpha_\pi(C^{\prime}),$ and 
$\alpha_\pi(C')\prec \alpha_\pi(C)$. On the other hand we have
$\alpha_\pi(C^{\prime})\in \mathrm{cl}\big(C'\setminus 
\alpha_\pi(C^{\prime})\big)\subseteq 
\mbox{cl}(C\setminus c)=\mathrm{cl}(C\setminus
\alpha_\pi(C)),$ so $\alpha_\pi(C)\prec \alpha_\pi(C^{\prime}),$ a contradiction.
\end{proof}
\begin{corollary}\label{puta} The set\,
$\mathcal{G}_{u}:=\big\{\partial e_{C}:\, C\in 
\mathfrak{C}(\M)\}$ is a minimal universal Gr\"obner basis of the ideal
\,$\Im_{\chi}(\M)$.
\end{corollary}
\begin{proof} From Theorem~\ref{red},
the reduced Gr\"obner bases constructed for the different orders $\prec$
are all contained in $\mathcal{G}_{u}$. We prove the minimality by
contradiction.  Let $C_0=\{i_1,\dotsc,i_m\}$
be a circuit of $\M$ and let $\pi\in \os_n$ be a permutation such that 
$\pi^{-1}(i_j)=j,$ $j=1,\dotsc, m$.
Then $\mathcal{G}^{\prime}_u:=\{\partial e_{C}:\, C\in 
\mathfrak{C}\setminus C_0\}$ 
it is not a Gr\"obner basis because
$\mathrm{lp}_{\prec_\pi}(\partial e_{C_0})=e_{C_0\setminus i_1}$ is not 
in 
$\mathrm{Lt}_{\prec_\pi}(\mathcal{G}^{\prime}_u)$.
\end{proof}
To finish this section we give an important characterization of the no 
broken circuit 
bases of the $\chi$-algebras in terms 
of the Gr\"obner bases of their ideals.
\begin{definition}
{\em
Consider a permutation $\pi \in \os_n$ and the associated 
re-ordering $<_\pi$ of $[n]$. 
When the $<_\pi$-smallest element $\mathrm{inf}_{<_\pi}(C)$  
of a circuit $C\in \mathfrak{C}(\M),$ $|C|>1,$ is
deleted, the remaining  set,  
$C\setminus \mathrm{inf}_{<_\pi}(C),$\, is called a
\emph{$\pi$-broken circuit} of $\M$. We say that
$$\mbox{$\pi$-$\boldsymbol{nbc}(\M)$}:=\{e_X:\, X\subseteq [n]\,\,
\mbox{contains no}\,\, \mbox{$\pi$-broken circuit of}\,\, \M\}$$ is  
the \emph {$\pi$-no broken circuit bases} of 
${\mathbb A}_{\chi}(\M)$. As the algebra 
${\mathbb A}_{\chi}(\M)$ does not depend
of the ordering of the elements of
$\M$ it is clear that $\pi$-$\boldsymbol{nbc}(\M)$ is a no broken circuit
bases of ${\mathbb A}_{\chi}(\M)$.
}
\end{definition}
\begin{corollary}\label{David3}Let $\mathcal{B}$ be
a basis  of the module  $ {\mathbb A}_{\chi}(\M) $. Then are equivalent:
\begin{namelist}{xxxxxxxxxx}
\item[~$(\ref{David3}.1)$]$\mathcal{B}$ is the canonical 
basis $\mathcal{B}_\prec,$
for some term order $\prec$ on\, $\mathbb{T}(\Af)$. 
\item[~$(\ref{David3}.2)$]$\mathcal{B}$ is the 
$\pi$-no broken circuit bases $\pi$-$\boldsymbol{nbc}(\M)$, for some
permutation $\pi\in \os_n$.
\item[~$(\ref{David3}.3)$]$\mathcal{B}$ is the canonical basis 
$\mathcal{B}_{\mathcal{G}_r},$ 
for some reduced Gr\"obner basis $\mathcal{G}_r$ of the  ideal 
$\Im_{\chi}(\M)$.
\end{namelist}
\end{corollary}
\begin{proof}\, $(\ref{David3}.1)\Rightarrow (\ref{David3}.2)$\,  
Let $\prec$ be a term order of $\mathbb{T}(\Af)$. Since
$\mathcal{G}_{u}$ is a universal Gr\"obner basis of
$\Im_{\chi}(\M)$ (see  Corollary \ref{puta}) it is trivially a 
Gr\"obner basis relatively to $\prec$.
We have already remarked that the leading term of $\partial e_C$ is 
$e_{C\setminus c}$ where $c=\mbox{inf}_{<_{\pi_\prec}}(C)$.
From Proposition~\ref{David} we conclude that 
$\mathcal{B}_{\prec}=\pi_\prec$-$\boldsymbol{nbc}(\M)$.
\par
\noindent{$(\ref{David3}.2)\Rightarrow (\ref{David3}.3)$} 
Suppose that  $\mathcal{B}=\pi$-$\boldsymbol{nbc}(\M)$. 
Let $\prec_\pi$ be the degree lexicographic order of
$\mathbb{T}$ determined by the permutation $\pi\in \os_n$. 
Note that $\pi_{\prec_\pi}=\pi$.
From Theorem~\ref{red} we know that $(\mathcal{G}_r)_
{\prec_\pi}=\big\{\partial e_{C}:\, C\in 
\mathfrak{C}_{\prec_\pi}\}$
is the  reduced Gr\"obner
basis of\, $\Im_{\chi}(\M)$ with respect to the term order $\prec_\pi$.
Then $\mathcal{B}$ is the canonical basis for the
reduced Gr\"obner basis $(\mathcal{G}_{r})_{\prec_\pi}$.
\par
\noindent{$(\ref{David3}.3)\Rightarrow (\ref{David3}.1)$}\,
It is a consequence of Proposition~$\ref{David}$ and
Remark~\ref{David2}.
\end{proof}
\section{Diagonal bases of $\chi$-algebras}
\begin{proposition}\label{con}
Let $\mathbb{A}_\chi({\M})$ be a $\chi$-algebra with the associated map 
$\chi: 2^{[n]}\to \K$. For any non 
loop element $x$ of
$\M([n]),$ we define the two maps:
\begin{equation}\label{eq1}
\chi_{\M \setminus x} : 2^{[n]\setminus x}\to \K\hspace{5mm}  
\mbox{by}\hspace{5mm}
\chi_{\M \setminus x}(X)=\chi(X)\hspace{5mm}  \mbox{and}
\end{equation} 
\begin{equation}\label{eq2}
\chi_{\M /x} : 2^{[n]\setminus x}\to \K \hspace{5mm}  
\mbox{by}\hspace{5mm}
\chi_{{\M /x}}(X)=\chi(X\circ x).
\end{equation} 
 There are two $\chi$-algebras,\,
$\mathbb{A}_{\chi_{\M /x}}({\M/x})$\, and\,
$\mathbb{A}_{\chi_{\M \setminus x}}({\M\setminus x})$,\, associated to the
maps   \,  $\chi_{\M \setminus x}$ and\,  $\chi_{\M /x}$, respectively.
\end{proposition}
\vspace{2mm}
\noindent{\em Proof}.\,
From   $(\ref{FL1}.1)$ we know that $\chi(X)\not =\emptyset$ iff
$X\in \mbox{IND}(\M)$.
The deletion case being trivial, we will just prove the contraction case.
We have to show that $\chi_{\M /x}$  verifies properties $(\ref{FL1}.1)$
and $(\ref{FL1}.2)$. The first property is verified since a set $I$ is
independent in 
$\M /x$  iff
$I\cup x$ is independent in $\M$. To see that the second property 
is also verified, let
$U$ and $U'$ be two unidependents sets of $\M /x$. I.e.,  iff $U\cup x$ and $U'\cup x$ 
are two unidependents sets of $\M$. From $(\ref{FL1}.1)$ we
know that $$\partial e_{U\cup x}=\varepsilon_{_{{U\cup x},{U'\cup x}}}
(\partial e_{{U'\cup x}})e_{U\setminus U'}\quad\mbox{where}\quad
\varepsilon_{_{{U\cup x},{U'\cup x}}}\in \K^*.$$ Let
$\partial^\prime$ be the $\chi_{\M /x}$-boundary, i.e., the linear mapping $\partial^\prime:
\Af/\langle e_x \rangle
\to
\Af/\langle e_x \rangle$ such that for ever $e_i\in E\setminus x$ we 
have $\partial^\prime e_i=1$,  $\partial^\prime e_\emptyset =1$
and for every monomial
$e_X,\, x\not \in X$ and
$X\not =
\emptyset$,
$$\partial^\prime e_{X}=\sum_{p=1}^{p=m}(-1)^p
\chi_{{\M /x}} (X\setminus i_p)e_{X\setminus i_p}=\sum_{p=1}^{p=m}(-1)^p
\chi(X\setminus i_p\circ x)e_{X\setminus i_p}.$$
To finish the proof
we will show
that there is a scalar $\tilde{\varepsilon}_{_{U,U'}}\in \K^*$ such that 
$$\partial^\prime e_U=\tilde{\varepsilon}_{_{U,U'}}
(\partial^\prime e_{U'})e_{U\setminus U'}.$$
Let $X,X'\subseteq [n]$ be two disjoint subsets. From
Definition~\ref{algebra} we known that  
$$e_X\cdot e_{X'}=
\beta_{X,X'}e_{X\cup
X'},\,\, \mbox{where}\,\,\, \beta_{_{X,X'}}=\prod 
\beta_{i,j},\,\, (e_i\in X,e_j\in X'\,\, 
\mbox{and}\,\,  i>j).$$ So we have
with\, $U=(i_1,\dotsc, i_{m})$ and $U'=(j_1,\dotsc, j_{k})$, $U\cap
U^\prime=\emptyset$, $x\not\in U\cup U^\prime$:
 $$\pm\partial e_{U \cup x}= \sum_{p=1}^{p=m}(-1)^p
\chi(U\setminus i_p\circ x)e_{U\cup x\setminus
i_p}+(-1)^{m+1}\chi(U)e_U,$$
$$\partial^\prime e_{U}= \sum_{p=1}^{p=m}(-1)^p
\chi(U\setminus i_p\circ x)e_{U\setminus i_p},$$
$$\pm(\partial e_{U' \cup x})e_{U\setminus U'}=
\sum_{p=1}^{p=k}(-1)^p
\chi(U'\setminus j_p\circ x)\cdot\beta\cdot e_{U\cup x\setminus
j_p}+(-1)^{k+1}\chi(U')\cdot\beta^\prime\cdot e_{U},$$
where\, $\beta=\beta_{_{U'\cup x\setminus j_p, U\setminus
U'}}$\, and\, $\beta^\prime=\beta_{_{U',U\setminus U'}}$.
 $$(\partial^\prime e_{U' })e_{U\setminus U'}=\sum_{p=1}^{p=k}(-1)^p
\chi(U'\setminus j_p\circ x)\cdot\beta_{_{U'\setminus j_p,U\setminus U'}}
\cdot e_{U\setminus j_p}.$$
After remarking that
$$\beta_{_{U'\cup x\setminus j_p,U\setminus U'}}
\beta_{_{U'\setminus j_p,U\setminus U'}}^{-1}=
\beta_{_{x,U\setminus U'}}$$  does not
depend on $j_p,$ we can deduce that 
$$\partial^\prime e_U=\tilde{\varepsilon}_{_{U,U'}}(\partial^\prime
e_{U'})e_{U\setminus U'}\quad\mbox{with}\quad \tilde{\varepsilon}_{_{U,U'}}=\pm\varepsilon_{_{{U\cup
x},{U'\cup x}}}\cdot{\beta_{x,U\setminus
U'}^{-1}}.$$
\qed
\vspace{2mm}
\begin{proposition}\label{l} For every    non loop element $x$ of
$\M([n]),$
there is a unique monomorphism of\, vector spaces,\,
$\mathfrak{i}_x:  \mathbb{A}(\M\setminus x)\to
\mathbb{A}(\M),$ such that such that for every
$I\in
\mathrm{IND}(\M\setminus x),$ we have\,
$\mathfrak{i}_x(e_I)=e_I$.  
  \end{proposition}
\begin{proof}
 By a 
reordering  of the elements of the matroid
$\M$ we can suppose
that
     $x=n$.
It is clear that 
$$\mbox{NBC}(\M\setminus x)=\big\{X:
     X\subseteq [n-1] 
\,\,\mbox{and}\,\, X\in \mbox{NBC}(\M)\big\},$$
so the proposition is a consequence of Equation~(\ref{eq1}).
\end{proof}
\begin{proposition} For every non loop element $x$ of
$\M([n]),$ 
there is a unique epimorphism of vector spaces, 
$\boldsymbol{\mathfrak{p}}_x: 
\mathbb{A}({\mathcal{M}})\to
\mathbb{A}({\mathcal{M}}/x),$
 such that, for every
$e_I,$ $I\in
\mathrm{IND}({\mathcal{M}}),$ we have
\begin{equation}\label{p_x}
\boldsymbol{\mathfrak{p}}_x(e_I):=
\begin{cases}
\vspace{2mm}
e_{I\setminus x} & 
   \mbox{\em if} \hspace{5mm} x\in I,\\
\frac{\chi ( I\setminus y , x)}{\chi (I\setminus y,y)}\,e_{I\setminus y} & 
   \mbox{\em if there is }  y\in I 
 \hspace{2mm}\mbox{\em  parallel 
to $x,$}\\
0 & \mbox{\em otherwise}. 
\end{cases}
\end{equation}
\end{proposition}
\noindent{\em Proof.}\, From Remark~\ref{rem1},  it is enough to
prove that $\mathfrak{p}_x({\partial}e_{U})=0,$ 
for all unidependent $U=(i_1,\dotsc,i_{m})$.
We recall that if $x\in U$ then
${{U}}\setminus x$ is a unidependent set of
${\mathcal{M}}/x$. 
There are only the  following four cases:
\begin{namelist}{xx}
\item[~~$\circ$] If $U$ contains $x$ but no $y$ parallel to $x$ then:
\begin{align*}
\pm\mathfrak{p}_x({\partial}e_{U})=\mathfrak{p}_x(
(-1)^{m}\chi(U\setminus x)e_{U\setminus x}+
\sum_{i_p\in U\setminus x}(-1)^p
\chi(U\setminus \{i_p,x\}\circ x)e_{U\setminus i_p}))\\
=\sum_{i_p\in U\setminus x}(-1)^p
\chi(U\setminus \{i_p,x\}\circ x)e_{U\setminus \{i_p,x\}}=0\hspace{30mm}
\end{align*}
from Proposition \ref{con}.
\item[~~$\circ$] If $U$ does not 
contain $x$ but contains a $y$ parallel to $x$ then:
\begin{align*}
\pm\mathfrak{p}_x({\partial}e_{U})=\mathfrak{p}_x\big(
(-1)^{m}\chi(U\setminus y)e_{U\setminus y}+
\sum_{i_p\in U\setminus y}(-1)^p
\chi(U\setminus \{i_p,y\}\circ y)e_{U\setminus i_p}\big)\\
=
\sum_{i_p\in U\setminus y}(-1)^p
\chi(U\setminus \{i_p,y\}\circ y)\frac{\chi(U\setminus \{i_p,x\}\circ x)}
{\chi(U\setminus \{i_p,y\}\circ y)}e_{U\setminus \{i_p,y\}}=
0\hspace{1mm}
\end{align*}
like previously since $U\setminus y$ is again a unidependent of $\M/x$.
\item[~~$\circ$] If $U$   contains $x$ and a $y$ parallel
  to $x$ then:
\begin{align*}
\pm\mathfrak{p}_x({\partial}e_{U})=\mathfrak{p}_x(
\chi(U\setminus\{x,y\}\circ y)e_{U\setminus
x}-\chi(U\setminus\{x,y\}\circ x) e_{U\setminus
y})\\
=
\chi(U\setminus \{x,y\}\circ y)\frac{\chi(U\setminus \{x,y\}\circ x)}
{\chi(U\setminus \{x,y\}\circ y)}e_{U\setminus 
\{x,y\}}-\chi(U\setminus\{x,y\}\circ x)e_{U\setminus  
\{x,y\}}=
0.
\end{align*}
 \item[~~$\circ$] If $U$ does not 
contain $x$ nor a $y$ parallel to $x$ then:
\begin{equation*}
\hspace{25mm}\mathfrak{p}_x({\partial}e_{U})=\mathfrak{p}_x\big(
\sum_{i_p\in U}(-1)^p
\chi(U\setminus i_p)e_{U\setminus i_p}\big)=0.
\hspace{25mm}\qed
\end{equation*}
\end{namelist}
 \begin{theorem}\label{thm: There is}
For every  element $x$ of a simple $\M([n]),$    there is a 
splitting short  exact
sequence of vector spaces
\begin{equation}\label{exact}
    0\to \mathbb{A}(\M\setminus
x)\stackrel{\mathfrak{i}_x}\longrightarrow 
    \mathbb{A}(\M)\stackrel{\mathfrak{p}_x}\longrightarrow 
    \mathbb{A}(\M/x)\to 0.
    \end{equation}
\end{theorem}
\begin{proof}
From the definitions we know that the composite map
$\mathfrak{p}_x {}^{_\circ}\mathfrak{i}_x$,  is the null map
 so $\mathrm{Im}(\mathfrak{i}_x)\subseteq \mathrm{Ker}(\mathfrak{p}_x)$. 
    We will prove the equality  $\dim(\mathrm{Ker}(\mathfrak{p}_n))=\dim
(\mathrm{Im}(\mathfrak{i}_n))$. By a reordering  of the elements of $[n]$
 we can suppose
that
     $x=n$.
The minimal 
 broken circuits of 
$\M/n$ are the minimal sets $X$ such that either $X$
or $X\cup \{n\}$ is a broken circuit of $\M$ (see the Proposition 3.2.e
of \cite{Bry}).  Then
\begin{equation*}
\mbox{NBC}(\M/n)=\big\{X: X\subseteq 
 [n-1]~~\mbox{and}~~X\cup\{n\}\in \mbox{NBC}(\M)
\big\}\,\,\,\,\,\,\mbox{and}
\end{equation*}
\begin{equation}\label{NBC}
\mbox{NBC}(\M)=
\mbox{NBC}(\M\setminus n)\biguplus
\big\{I\cup
n: I\in
\mbox{NBC}(\M/n)\big\}.
\end{equation}
So $\dim(\mathrm{Ker}(\mathfrak{p}_n))=\dim
(\mathrm{Im}(\mathfrak{i}_n))$. There is a  
 morphism of modules 
$$\mathfrak{p}^{-1}_n: \mathbb{A}(\M/n)\to
\mathbb{A},~~\mbox{where}~~~\mathfrak{p}^{-1}_n([I]_{\mathbb{A}(\M/n)})
:=[I\,\cup
\,n]_{\mathbb{A}},\, \forall I\in
\mathrm{NBC}(\M/n).$$
It is clear that the composite map $\mathfrak{p}_n{}^{_\circ}
\mathfrak{p}^{-1}_{n}$ is the identity map. From Equation~(\ref{NBC}) we
conclude that the exact sequence~(\ref{exact}) splits.
\end{proof}
Similarly to  \cite{Sz} (see also \cite{BV}), we now construct, making use  
of iterated contractions,  the dual bases 
 $\boldsymbol{nbc}_{\,\ell}^*=(b^*_i)$ of the bases  
$\boldsymbol{nbc}_{\,\ell}:=(b_j)$ of the vector space $\mathbb{A}_\ell$. 
 More precisely
$\boldsymbol{nbc}^*_{\,\ell}$ is the basis of $\mathbb{A}_\ell^*$ 
the vector space of the linear forms such that 
$\langle b^*_i, b_j\rangle=\delta_{ij}$ (the Kronecker delta).\par
We associate to  the ordered independent set 
$I^\sigma:=(i_{\sigma(1)},\dotsc,  i_{\sigma(p)})$ of\,
${\mathcal{M}}$ the linear form on 
$\mathbb{A}_\ell,$ 
$\boldsymbol{\mathfrak{p}}_{I^\sigma}:
\mathbb{A}_\ell\to  \K$, defined as the composite of the maps\,
$\boldsymbol{\mathfrak{p}}_{e_{i_{
\sigma(p)}}}$, $\boldsymbol{\mathfrak{p}}_{e_{i_{
\sigma(p-1)}}}, \dotsc,
\boldsymbol{\mathfrak{p}}_{e_{i_{\sigma(1)}}}$, i.e.,
\begin{equation}\label{resid}
\boldsymbol{\mathfrak{p}}_{I^\sigma}:=
\boldsymbol{\mathfrak{p}}_{e_{i_{
\sigma(1)}}}\!{}^{_\circ}\boldsymbol{\mathfrak{p}}_{e_{i_{\sigma(2)}}}\!
{}^{_\circ}
\cdots\,{}^{_\circ}
\boldsymbol{\mathfrak{p}}_{e_{i_{\sigma(p)}}}.
\end{equation}
   We call
$\boldsymbol{\mathfrak{p}}_{I^\sigma}$  the  \emph{iterated residue} with 
respect to the ordered independent
set $I^\sigma$. We remark that the map 
$\boldsymbol{\mathfrak{p}}_{I^\sigma}$ depends on the 
order  chosen on $I^\sigma$ and not only on
the  underlying set  $I$. We
associate to $I^\sigma$ the flag of flats of
${\mathcal{M}},$
$$
\textbf{Flag}(I^\sigma):=\cl\big(\{i_{\sigma(p)}\}\big) \subsetneq 
\cl\big(\{i_{\sigma(p)}, i_{\sigma(p-1)}\}\big)
\subsetneq 
\cdots   \subsetneq  \cl\big(\{i_{\sigma(p)}, \dotsc, i_{\sigma(1)}\}\big).
$$
\begin{proposition}\label{resi}
Let $J\in 
\mathrm{IND}_\ell({\mathcal{M}})$ then we have\, 
$\boldsymbol{\mathfrak{p}}_{I^\sigma}(e_{J})\not=0$ iff there is a unique 
permutation $\tau \in \os_\ell$  such that
$\textbf{Flag}\,(J^\tau) =
\textbf{Flag}\,(I^\sigma)$.  And in this case we have 
$\boldsymbol{\mathfrak{p}}_{I^\sigma}(e_{J})=
\chi (I^\sigma )/\chi (J^\tau)$. In particular we have\,
$\boldsymbol{\mathfrak{p}}_{I^\sigma}(e_{I})=1$ for any independent set 
$I$ 
and any permutation $\sigma$.
\end{proposition}
\begin{proof}
The first equivalence is  easy to prove in both direction.
To obtain the expression of $\boldsymbol{\mathfrak{p}}_{I^\sigma}(e_{J})$
we just need to iterate $\ell$ times the residue. This gives:
\begin{align*}
\boldsymbol{\mathfrak{p}}_{I^\sigma}(e_{J})=
\frac{\chi ( J\setminus j_{\tau (\ell)} \circ i_{\sigma(\ell)})}
{\chi ( J\setminus j_{\tau (\ell)}\circ j_{\tau (\ell)} )}
\times \frac{\chi ( J\setminus \{j_{\tau (\ell)}, j_{\tau (\ell
-1)}\}\circ   i_{\sigma (\ell -1)}\circ i_{\sigma(\ell)})}
{\chi ( J\setminus \{j_{\tau (\ell)}, j_{\tau (\ell -1)}\}\circ
j_{\tau (\ell-1)}\circ i_{\sigma(\ell )})} \times \dotsb \\
\dotsb \times
\frac{\chi (I^\sigma)}{\chi (j_{\tau (1)}\circ I^\sigma \setminus 
i_{\sigma(1)})}.
\end{align*}
After simplification we obtain the announced formula. The
last result is clear.
\end{proof}
\begin{remark}\label{rem}
{\em
The fact that $\boldsymbol {\mathfrak{p}}_{I^\sigma}(e_{J})$ is null 
depends on the permutation $\sigma$. For example, for any
simple matroid of rank 2  we have  
$\boldsymbol {\mathfrak{p}}_{13}(e_{12}) =0$ and 
$\boldsymbol{\mathfrak{p}}_{31}(e_{12})\not =0$.  But if $\boldsymbol
{\mathfrak{p}}_{I^\sigma}(e_{J})
\not=0$ then its value does not depend on $\sigma$. 
We mean by this that if there 
are two permutations $\sigma $ and $\sigma '$ such that 
$\boldsymbol{\mathfrak{p}}_{I^\sigma}(e_{J})\not =0$
and $\boldsymbol{\mathfrak{p}}_{I^{\sigma '}}(e_{J})\not =0$ then 
$\boldsymbol{\mathfrak{p}}_{I^\sigma}(e_{J})=
\boldsymbol{\mathfrak{p}}_{I^{\sigma
'}}(e_{J})$.
}
\end{remark}
\begin{definition}[\cite{Sz}]\label{diagonal}
 {\em
We
say that the subset $\mathbb{I}_\ell\subseteq \big\{[I]_{\mathbb{A}}: I\in 
\mbox{IND}_\ell({\mathcal{M}})\}$ is
 a
\emph{diagonal basis} of   $\mathbb{A}_\ell$ if and only if the following
three conditions hold:
\begin{namelist}{xxxxxx}
\item[~$(\ref{diagonal}.1)$] For every $[I]_{\mathbb{A}}\in 
\mathbb{I}_\ell$ there is a fixed
permutation of the set $I$ denoted $\sigma_I\in \os_\ell;$
 \item[~$(\ref{diagonal}.2)$] $\big|\mathbb{I}_\ell| 
\geq\mbox{dim}(\mathbb{A}_\ell);$
\item[~$(\ref{diagonal}.3)$]  For every $[I]_{\mathbb{A}}, 
[J]_{\mathbb{A}}\in
\mathbb{I}_\ell$ and every permutation
$\tau
\in \os_\ell,$ the equality
$\textbf{Flag}\,(J^\tau) =
\textbf{Flag}\,(I^{\sigma_I})$ implies $J= I$.
\end{namelist}
}
\end{definition}
\begin{theorem}\label{dbases} Suppose that $\mathbb{I}_\ell$ is a
diagonal  basis  of   $\mathbb{A}_\ell$.
Then $\mathbb{I}_\ell$ is a basis
of  $\mathbb{A}_\ell$ and
\,$\mathbb{I}_\ell^*:=\{\boldsymbol{\mathfrak{p}}_{I^{\sigma_I}}: 
[I]_\mathbb{A}\in \mathbb{I}_\ell\}$ is the dual 
 basis  of\, $\mathbb{I}_\ell$. 
\end{theorem}
\begin{proof} 
Pick two elements $[I]_\mathbb{A}, [J]_\mathbb{A}
\in\mathbb{I}_\ell$. Note that
$\mathfrak{p}_{I^{\sigma_I}}(e_{J})=\delta_{IJ}$ 
(the Kronecker delta), from Condition~$(\ref{diagonal}.2)$ and
Proposition~\ref{resi}.
 The elements of $\mathbb{I}_\ell$ are linearly independent: suppose that
$[J]=\sum\zeta_j [I_j],\,  \zeta_j\in \K\setminus 0;$ then
$1=\boldsymbol{\mathfrak{p}}_{J^{\sigma_J}}([J])=
\boldsymbol{\mathfrak{p}}_{J^{\sigma_J}}\big(\sum\zeta_j
[I_j]\big)=0,$ a contradiction. It is clear also that $\mathbb{I}_\ell^*$ 
is the dual 
 basis  of\, $\mathbb{I}_\ell$. 
\end{proof}
The following result gives an interesting explanation 
of results of \cite{cor, cor2}.
\begin{corollary}\label{that}
$\boldsymbol{nbc}_\ell({\mathcal{M}})$  is a diagonal basis of
$\mathbb{A}_\ell$ where $\sigma_I$ is the identity  for every 
$[I]_{\mathbb{A}}
\in
\boldsymbol{nbc}_\ell({\mathcal{M}})$.
  For a given $[J]_{\mathbb{A}}\in \mathbb{A}_\ell,$
suppose
that
\begin{namelist}{xxxxxxxx}
\vspace{1mm}
\item[$(\ref{that}.2)$]$[J]_{\mathbb{A}}
=\sum
\xi (I,J) [{I}]_{\mathbb{A}},$ where
\,\,$[I]_{\mathbb{A}}\in\boldsymbol{nbc}_\ell({\mathcal{M}})$\, and 
\,\, 
$\xi (I,J) \in \K$.
\vspace{1mm}
\end{namelist}
Then  are equivalent:
\begin{namelist}{xxxxx}
\item[~~$\circ$]$\xi (I,J)\not =0,$
\item[~~$\circ$]$\textbf{Flag}\,(I) =
\textbf{Flag}\,({J}^{\tau})$ for some permutation $\tau$.
\end{namelist} 
If $\xi (I,J)\not =0$ we have $\xi (I,J)=\frac{\chi (I)}
{\chi (J^\tau)}$. In particular
if $\mathbb{A}$  is the Orlik-Solomon algebra then
$\xi(I,J)=\sgn(\tau)$.
\end{corollary}
\begin{proof}
By hypothesis $(\ref{diagonal}.1)$ and $(\ref{diagonal}.2)$ are true.
 We claim that 
$\boldsymbol{nbc}_{\,\ell}({\mathcal{M}})$ verifies
$(\ref{diagonal}.3)$. Suppose for a contradiction that $J\not = I,$
$[J]_\mathbb{A}, [I]_\mathbb{A}\in 
\boldsymbol{nbc}_{\,\ell}({\mathcal{M}})$ and there is $\tau \in 
\os_\ell,$ such that $\textbf{Flag}\,(J^\tau) =
\textbf{Flag}\,(I)$. Set $I=(i_1,\dotsc,  i_\ell)$ and  
$J=(j_{\tau(1)},\dotsc,  j_{\tau({\ell})}),$  
and suppose that $j_{\tau(m+1)}=i_{m+1}, 
\dotsc, j_{\tau(\ell)}=i_\ell$ and $i_m\not =j_{\tau({m})}$. 
Then there is a circuit
$C$ of
$\M$ such that $$i_m, j_{\tau(m)}\in 
C\subseteq \{i_m, j_{\tau(m)}, i_{m+1}, i_{m+2}, \dotsc
,i_\ell\}.$$
If $j_{\tau(m)}<i_m$ [resp. $i_m<j_{\tau(m)}$] we conclude that $$I\not 
\in 
\mbox{NBC}_\ell({\mathcal{M}})\quad\mbox{
[resp.}\quad J\not \in \mbox{NBC}_\ell({\mathcal{M}})] $$ a 
contradiction. So $\boldsymbol{nbc}_{\,\ell}({\mathcal{M}})$  
is a diagonal basis of
$\mathbb{A}_\ell$.\par
From Theorem~\ref{dbases} we 
conclude that   
$\boldsymbol{nbc}_\ell^*:=\big\{\boldsymbol{\mathfrak{p}}_{I}: 
[I]_\mathbb{A}\in
\boldsymbol{nbc}\}$ is the dual 
 basis  of\, $\boldsymbol{nbc}$. Suppose now that 
$[J]_{\mathbb{A}}=\sum{\xi}_{I}[{I}]_{\mathbb{A}},$ where
$[I]_{\mathbb{A}}\in\boldsymbol{nbc}_\ell({\mathcal{M}})$ and 
${\xi}_I\in
k$. Then ${\xi}_{I}=\boldsymbol{\mathfrak{p}}_{I}(e_{J})$
and the remaining follows from Proposition~\ref{resi}.
\end{proof}
Making full use of the matroidal notion of iterated residue, 
see Equation~(\ref{resid}),
we are able to prove the following result very close to Proposition~ 2.1 
of
\cite{Sz2}.
\begin{proposition}\label{machain}
Consider the set  of  vectors  $\mathcal{V}:=\{v_1,\dotsc,  v_k\}$ 
in the plane $x_d=1$ of\, $\K^d$. Set $\A_\K:=\{H_i: 
H_i=\Ker(v_i)\subseteq (\K^d)^*,\,
i=1,\dotsc, k\}$  and let\, $\OT(\A_\K)$ be its Orlik-Solomon-Terao 
corresponding
algebra.  Fix a diagonal basis  
$\mathbb{I}_\ell\subseteq \{[I]_{\mathbb{A}}: I\in \mathrm{IND}_\ell(\M)\}
$
of  $\mathbb{A}_\ell$ and let\,
$\mathbb{I}^*_\ell=\{{\mathfrak{p}}_{I^{\sigma_I}}: 
[I]_{\mathbb{A}}\in \mathbb{I}_\ell\}
$
be the corresponding dual  basis. 
Then, for any\,  $e_J\in \mathbb{A}_\ell \setminus 0,$ 
we have 
$$\sum _{I\in \mathbb{I}_\ell} \boldsymbol 
{\mathfrak{p}}_{I^{\sigma_I}}(e_J)=\sum _{I\in
\mathbb{I}_\ell}
\big\langle {\mathfrak{p}}_{I^{\sigma_I}},e_J\big\rangle =1.$$
\end{proposition}
\begin{proof}
We have for any $\ell+1$-subset of $\mathcal{V},$
$\sum_{p=1}^{p=\ell+1}(-1)^p
\chi(U\setminus i_p)=0$. (This is the development of a 
determinant with two lines of 
1.)
For any rank $\ell$ unidependent $U=\{i_1,\dotsc,  
i_{\ell+1}\}$  of the matroid
$\mathcal{M} (\A_\K),$ we have 
$$\partial e_U=\sum_{p=1}^{p=\ell+1}(-1)^p
\chi(U\setminus i_p)e_{U\setminus i_p}.$$
Since  the sum of the coefficients in these relations is 
0 and that these relations are
generating, see Remark~\ref{rem1}, we can deduce that the 
sum of the coefficients in any
relation in
$\OT(\A_\K)$ is also equal to 0 which concludes the proof.
\end{proof}
\section{Examples}
In this section we will show on a small example
the different results of the three previous sections.\par
Consider  the the set of 6 points $\{p_1,\dotsc,  p_6\}$ in the
affine plane $z=1$ of
three dimensional real vector space $\R^3,$  whose 
coordinates are indicated in Figure~1. Set
$v_i:=\overrightarrow{(0,p_i)},\ i=1,\dotsc, 6$. And let
$\A$ be the corresponding hyperplane arrangement  of
$(\R^3)^*,$ $\A:=\{H_i=\Ker(v_i),\, i=1,\dotsc, 6\}$.
Let $\M(\A)$ [resp. $\boldsymbol{\M}(\A)$] be the corresponding rank three
[resp. oriented]  matroid. So like in Example \ref{chi}, the arrangement
$\A$ defines the three classical Orlik-Solomon type algebras:
the original Orlik-Solomon algebra $\OS (\M(\A))$ through $\M(\A)$, the 
Orlik-Solomon-Terao algebra $\OT(\A)$ directly from the $v_i$ and the 
Cordovil algebra  $\mathbb{A}(\boldsymbol{\M}(\A))$ from 
$\boldsymbol{\M}(\A)$.
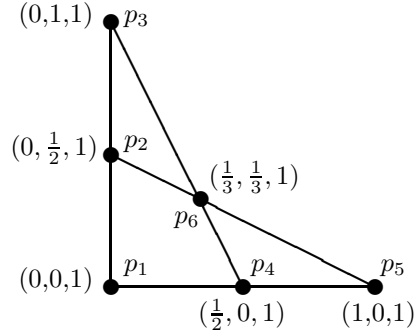
\begin{figure}
\begin{center}
\thicklines
\begin{picture}(-76,120)(100,10)
\put(15,16){$p_1$}
\put(15,110){$p_3$}
\put(15,63){$p_2$}
\put(63,16){$p_4$}
\put(112,16){$p_5$}
\put(34,34){$p_6$}
\put(10,10) {\circle*{6}}
\put(60,10) {\circle*{6}}
\put(10,60) {\circle*{6}}
\put(10,110) {\circle*{6}}
\put(44,43) {\circle*{6}}
\put(110,10) {\circle*{6}}
\put(-25,10){(0,0,1)}
\put(-28,60){$(0,\frac {1}{2},1)$}
\put(-25,110){(0,1,1)}
\put(43,-3){$(\frac {1}{2},0,1)$}
\put(97,-3){(1,0,1)}
\put(47,49){$(\frac {1}{3},\frac{1}{3},1)$}
\put(10,10) {\line(2,0){100}}
\put(10,10) {\line(0,1){100}}
\put(10,60) {\line(2,-1){100}}
\put(10,110) {\line(1,-2){50}}
\end{picture}
\vspace{3mm}
\caption{The rank 3 matroid on the set $\{p_{1},\dotsc, p_{6}\}$.}
\end{center}
\end{figure}
Let $\mathbb{A}_\chi$ be a  $\chi$-algebra on $\M(\A)$.
We know that $$\boldsymbol{nbc}_{\,3}=\{e_{124},e_{125},e_{126},
e_{134},e_{135},e_{136}\}$$ together with $\sigma _{124}=\sigma _{125}=
\sigma _{134}=\sigma _{135}=\sigma _{136}=\sigma _{156}=
\mbox{id}$ is a diagonal basis of $\mathbb{A}_3,$ from
Corollary~\ref{that}. Directly from  the
Definition~\ref{diagonal} we see that  $\B_3=\{ e_{124},e_{125},e_{134},
e_{135},e_{136},e_{156}\}$ with
$\sigma _{124}=\sigma _{134}=\sigma _{135}=\sigma _{136}=\sigma _{156}=
\mbox{id}$ and $\sigma _{125}=(132)$ is also a diagonal basis of 
$\mathbb{A}_3$. We will look at expressions 
on the  basis $\boldsymbol{nbc}_3$ (resp. $\B_3$) of the vector space 
$\mathbb{A}_3,$
of some
elements of the type $e_B,$ $B$ basis of $\M(\A),$   for the three 
$\chi$-algebras  of  Example~\ref{chi}. Especially, we will verify   as
 stated in Remark~\ref{rem}\, that
 $\boldsymbol{\mathfrak{p}}_{125^{\id}}(e_{235})=
\boldsymbol{\mathfrak{p}}_{125^{(132)}}(e_{235})$. Let also point out
that for
 the Orlik-Solomon-Terao algebra, we will have have $\sum _{I\in
 \B}\boldsymbol{\mathfrak{p}}_{I^{\sigma}}(e_{J})=1$
as proved in Proposition~\ref{machain}.
Finally recall that $\mathbb{T}$ is set of the monomials  of $\Af$  and set  $\mathbb{T}_\ell:=\{e_X\in \mathbb{T}:\, |X|=\ell\}$.
\begin{namelist}{xxxx}
\item[~$(a)$]
Let us first take the Orlik-Solomon algebra $\OS (\M(\A)):$\par
From Remark \ref{rem1} , the basis of $\OS (\M(\A))$ is simply the 
$\boldsymbol{nbc}$-bases:
$$\boldsymbol{nbc}(\M)=\mathbb{T}_0\cup\mathbb{T}_1\cup 
\boldsymbol{nbc} _2\cup \boldsymbol{nbc} _3,$$ with\,
$\boldsymbol{nbc} _2=\{e_{12}, e_{13}, e_{14}, e_{15}, e_{16},
e_{24}, e_{25}, e_{26}, e_{34}, e_{35}, e_{36}\}$,\,  and
$$\boldsymbol{nbc} 
_3=\{e_{124},e_{125},e_{126},e_{134},e_{135},e_{136}\}.$$
The basis of $\Im_{\chi} (\M(\A))$ is the union of the dependents and
of the boundaries of the 
inactive unidependents:
$$ \partial \boldsymbol{ uni} _3
\cup\boldsymbol{dep} _3\cup\partial \boldsymbol{uni} _4
\cup\mathbb{T}_4\cup\mathbb{T}_5\cup\mathbb{T}_6 $$ where 
$\partial \boldsymbol{uni} _3=
\{\partial e_{123},\partial e_{145},\partial e_{256},
\partial e_{346}\},$ 
$\boldsymbol{dep} _3=\{e_{123},
e_{145},e_{256},e_{346}\} $ and
$\partial \boldsymbol{uni} _4$ is the set
$$\{
\partial e_{1234},\partial e_{1235},\partial e_{1236}
,\partial e_{1245},\partial e_{1246},\partial e_{1256},\partial e_{1345}
,\partial e_{1346},\partial e_{1356},\partial e_{1456}
\}.
$$
Note that we have  
$$|\boldsymbol{nbc} _2|+|\partial \boldsymbol{uni} _3|=11+4=15=
\dim (\Af _2)$$ and $$|\boldsymbol{nbc} _3|+|\partial  \boldsymbol{uni}
_4| +|\boldsymbol{dep} _3|=6+10+4=20=
\dim (\Af _3).$$
Take first on $[n]$ the natural order. We have then for the leading
term ideal
$$\mbox{Lt}_<(\mathcal{G})=\langle e_{BC}:BC\,\, \mbox{broken
circuit}\rangle.$$ We obtain explicitly:
$$\mbox{Lt}_< (\mathcal{G})=\langle  e_{23},e_{45},e_{56},e_{46},
e_{246},e_{345},e_{356}\rangle.$$
Always for the natural order, from Theorem \ref{red}, we obtain 
for the reduced Gr\"obner basis:
$$\mathcal{G}_{r} =\big\{\partial e_{123},\partial e_{145},
\partial e_{256},\partial e_{346}\big\}.$$
If we take now  the term order\,
$\prec_\pi$\, on $\mathbb{T}(\Af)$,  defined by the permutation \,$\pi:=(234561),$  we get now:
$$\mbox{Lt}_\prec (\mathcal{G})=\langle e_{13},e_{15},e_{56},e_{46},
e_{146},e_{345},e_{165}\rangle,$$
and then for the corresponding reduced Gr\"obner basis:
$$\mathcal{G}_{r} =\big\{\partial e_{123},\partial e_{145},
\partial e_{256},\partial e_{346},\partial e_{2345}\big\}.$$ 
Finally from Corollary \ref{puta}, we get the minimal 
universal Gr\"obner basis 
$$\mathcal{G}_{u}=\big\{\partial e_{C}:\, C\in 
\mathfrak{C}(\M)\}.$$
We obtain explicitly:
$$\mathcal{G}_{u}=\{\partial e_{123},\partial e_{145},
\partial e_{256},\partial e_{346},\partial e_{1246},\partial e_{1356},
\partial e_{2345}
\}.$$
Now we will use the results of Section 4 to express pure elements
in different diagonal bases. Consider the  diagonal  basis 
$\boldsymbol{nbc}_3$ of the $\mathbb{K}$-vector space
$\OS (\M(\A))_3$. So
 we have:
$$e_{156}=\mbox{sgn}(165)e_{125}+ \mbox{sgn}(156)e_{126}=
-e_{125}+e_{126}$$
and
$$e_{235}=\sgn(325)e_{125} + \sgn(235)e_{135}=
-e_{125} + e_{135}.$$
For the diagonal basis $\B_3$ of the $\mathbb{K}$-vector space
$\OS (\M(\A))_3,$ 
 we have:
$$e_{126}=\mbox{sgn}(162)\mbox{sgn}(152)e_{125}+ \mbox{sgn}(126)e_{156}=
e_{125}+e_{156}$$
and
$$e_{235}=\mbox{sgn}(152)\mbox{sgn}(352)e_{125}+ 
\mbox{sgn}(235)e_{135}=
-e_{125}+e_{135}.$$
\item[~$(b)$] Let us take the Orlik-Solomon-Terao algebra $\OT(\A):$\par
For the different bases and Gr\"obner bases we obtain 
formally the same results.
There is in fact differences which are hidden by the operator $\partial$
(indeed  $\partial$ is function of $\chi$).\par
For the  diagonal  basis $\boldsymbol{nbc}_3$ 
of the $\K$-vector space
$\OT(\A)_3$  we have:
$$e_{156}=\frac {\det (125)}{\det (165)}e_{125}+ 
\frac {\det (126)}{\det (156)}e_{126}
=\frac {3}{2}e_{125}-\frac {1}{2}e_{126}$$
and
$$e_{235}=\sgn(325)e_{125} + \sgn(235)e_{135}=
-e_{125} + e_{135}$$
For the diagonal basis $\B_3$ of the $\mathbb{K}$-vector space
$\OT(\A)_3$ 
 we have:
$$e_{126}=\frac {\det (152)}{\det (162)}e_{125}+ 
\frac {\det (156)}{\det (126)}e_{156}=3e_{125}-2e_{156}.$$
and
$$e_{235}=\frac {\det (152)}{\det (352)}e_{125}+ 
\frac {\det (135)}{\det (235)}e_{135}=-e_{125}+2e_{135}.$$
\item[~$(c)$] 
Let us take the Cordovil $\Z$-algebra 
$\mathbb{A}(\boldsymbol{\M}(\A)):$\par
For the  diagonal  basis $\boldsymbol{nbc}_3$ 
of the $\mathbb{K}$-vector space
$\mathbb{A}(\boldsymbol{\M}(\A))_3$ 
 we have:
$$e_{156}=\chi (125)\chi (165)e_{125}+ 
\chi (126)\chi (156)e_{126}=e_{125}-e_{126}$$
and
$$e_{235}=\sgn(325)e_{125} + \sgn(235)e_{135}=
-e_{125} + e_{135}.$$
For the  diagonal basis $\B_3$ of the $\K$-vector space
$\mathbb{A}(\boldsymbol{\M}(\A))_3$
 we have:
$$e_{126}=\chi (152)\chi (162)e_{125}+ 
\chi (156)\chi (126)e_{156}=e_{125}-e_{156}
$$
and
$$e_{235}=\mbox{sgn}(152)\mbox{sgn}(352)e_{125}+ 
\mbox{sgn}(235)e_{135}=
-e_{125}+e_{135}.$$
\end{namelist}

\end{document}